\documentclass[10pt,twocolumn,letterpaper, capitalize]{article}
\usepackage{tikz}
\usepackage[review]{cvpr}      
\usepackage{float}
\definecolor{cvprblue}{rgb}{0.21,0.49,0.74}

\title{
    Blackout DIFUSCO:\\
    Exploring Continuous-Time Dynamics\\
    in Discrete Optimization}

\begin{document}
\maketitle
\begin{figure*}[t] 
    \centering
    \includegraphics[width=\textwidth,keepaspectratio]{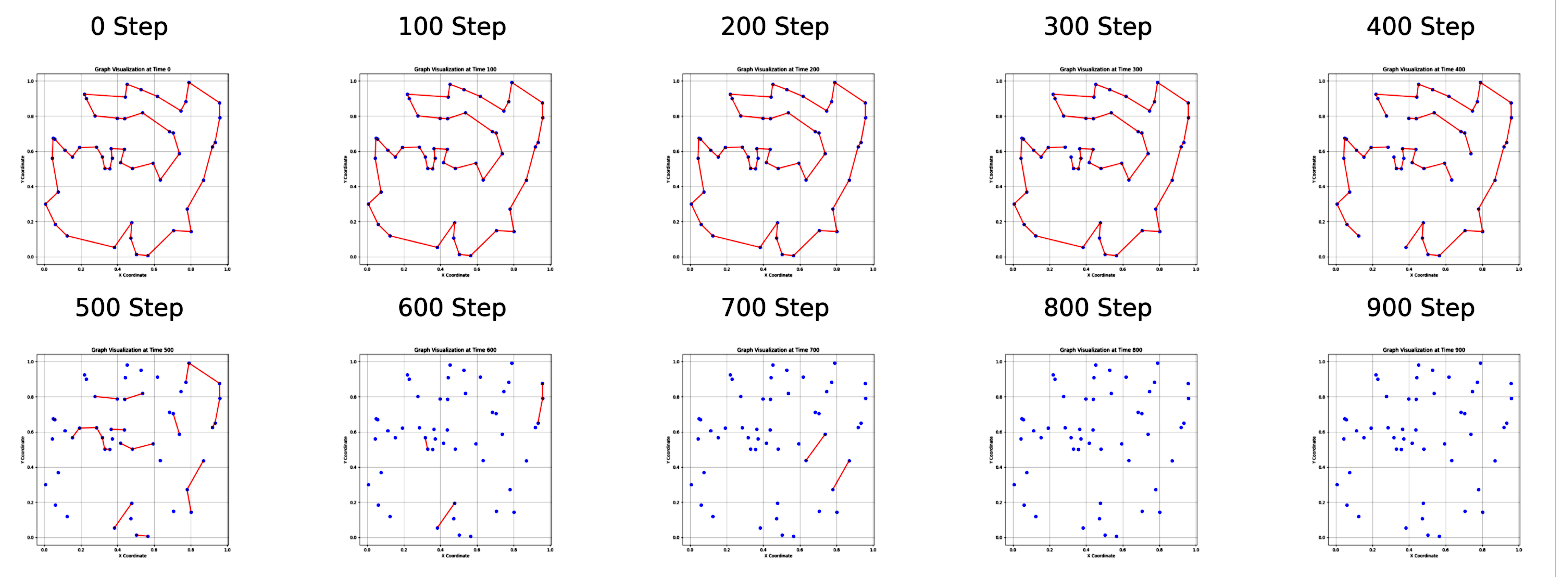}
    \caption{Visualization of the forward process in Blackout DIFUSCO. The images depict the progressive corruption of the adjacency matrix as the diffusion process moves from the initial state (\(x_0\)) to a fully corrupted state (\(x_T\)). Each frame corresponds to a specific timestep in the forward diffusion process.}
    \label{fig:forward_process}
\end{figure*}

\section*{Abstract}
This study explores the integration of Blackout Diffusion into the DIFUSCO framework for combinatorial optimization, specifically targeting the Traveling Salesman Problem (TSP). Inspired by the success of discrete-time diffusion models (D3PM) in maintaining structural integrity, we extend the paradigm to a continuous-time framework, leveraging the unique properties of Blackout Diffusion. Continuous-time modeling introduces smoother transitions and refined control, hypothesizing enhanced solution quality over traditional discrete methods. We propose three key improvements to enhance the diffusion process. First, we transition from a discrete-time-based model to a continuous-time framework, providing a more refined and flexible formulation. Second, we refine the observation time scheduling to ensure a smooth and linear transformation throughout the diffusion process, allowing for a more natural progression of states. Finally, building upon the second improvement, we further enhance the reverse process by introducing finer time slices in regions that are particularly challenging for the model, thereby improving accuracy and stability in the reconstruction phase. Although the experimental results did not exceed the baseline performance, they demonstrate the effectiveness of these methods in balancing simplicity and complexity, offering new insights into diffusion-based combinatorial optimization. This work represents the first application of Blackout Diffusion to combinatorial optimization, providing a foundation for further advancements in this domain. * The code is available for review at https://github.com/Giventicket/BlackoutDIFUSCO.
\section{Introduction} \label{sec:intro}

Combinatorial optimization (CO) involves finding the optimal solution from a finite, discrete solution space. These NP-hard problems, such as the Traveling Salesman Problem (TSP), Knapsack Problem, and Job Scheduling, require significant computational resources to solve exactly, motivating the development of efficient approximation methods.

Among CO problems, TSP is a widely studied example. It entails finding the shortest route to visit all cities exactly once and return to the starting point. The problem's factorial complexity makes brute-force solutions impractical, especially for large instances, necessitating heuristic and learning-based approaches.

To address NP-hard CO problems, supervised learning (SL), reinforcement learning (RL), and semi-supervised learning (SSL) methods have been explored. SL approaches generate approximate solutions using historical optimal data, often producing stable results via heatmap-based decoding \cite{joshi2019tsp, sun2023difusco, li2023t2t}. RL methods train agents to construct solutions sequentially or iteratively improve solutions through metaheuristics \cite{kool2019attention, bresson2021transformer, zheng2021rl-lkh}. However, RL approaches can suffer from inefficiencies during training due to large initial Optimality Gaps. Hybrid methods that integrate SL into RL training have been proposed to address this challenge.

Recently, DIFUSCO \cite{sun2023difusco} introduced diffusion generative models for CO problems, achieving state-of-the-art performance on TSP. DIFUSCO generates heatmaps representing edge probabilities and decodes them into feasible solutions. Building on this, the T2T model \cite{li2023t2t} introduced objective-guided sampling, further improving performance and efficiency while retaining DIFUSCO's non-autoregressive structure. DIFUSCO and T2T operate in a discrete time-space framework, with categorical-state diffusion (D3PM) outperforming continuous-state diffusion (DDPM) in discrete optimization tasks.

Inspired by these advancements, we explore Blackout Diffusion \cite{santos2023blackout}, a framework that bridges discrete state spaces and continuous-time modeling, for CO problems. Continuous-time diffusion, as demonstrated in "Score-Based Generative Modeling through Stochastic Differential Equations" \cite{song2021sde}, has shown performance gains in generative tasks compared to discrete-time methods like DDPM \cite{ho2020ddpm}. We hypothesize that applying continuous-time discrete-state diffusion to CO tasks could improve solution quality by enabling smoother state transitions.

This work represents the first attempt to combine DIFUSCO with Blackout Diffusion. By incorporating insights from diffusion scheduling and optimal observation time design, we aim to refine sampling strategy and solution quality. Although the performance was not satisfactory, our study lays a foundation for future research into continuous-time diffusion models for combinatorial optimization.

\subsection*{Contributions}

\begin{itemize} \item \textbf{Integration of DIFUSCO and Blackout Diffusion for Combinatorial Optimization}
\item \textbf{Sampling Optimization for Efficient Blackout Diffusion}
\item \textbf{Design of Heuristic Diffusion Processes}
\item \textbf{Exploration of Diffusion Models for Discrete State-Space Optimization}
\end{itemize}
\section{Related Work}
\label{sec:RelatedWork}

\subsection{TSP Formulation}

In the case of the Traveling Salesman Problem (TSP), the goal is to compute the shortest Hamiltonian cycle that visits all $N$ nodes in a Euclidean space exactly once and returns to the starting node. Here, $S$ is the set of all feasible tours, and $x \in S$ represents a specific tour. The objective function $f(x)$ is the total length of the tour, defined as:
\[
f(x) = \sum_{(i, j) \in x} d(i, j),
\]
where $d(i, j)$ is the Euclidean distance between nodes $i$ and $j$, and the summation is over all edges $(i, j)$ in the tour $x$.

It is important to note that $f(x)$ is only defined over the feasible set $S$, meaning $x$ must satisfy the constraints of forming a valid Hamiltonian cycle. In TSP, this ensures that every node is visited exactly once and the tour is closed.

\subsection{DIFUSCO: Graph-Based Diffusion Solvers for Combinatorial Optimization}

\paragraph{Diffusion Framework and Heatmap Representation.}  
DIFUSCO\cite{sun2023difusco} leverages diffusion models to solve combinatorial optimization problems by progressively corrupting and reconstructing adjacency matrices. The ground truth adjacency matrix, denoted as \(x_0\), encodes the presence or absence of edges in a binary or relaxed format. Through a forward diffusion process, \(x_0\) is gradually transformed into a fully corrupted state \(x_T\) over \(T\) timesteps. The reverse diffusion process then reconstructs the structured adjacency matrix by iteratively refining the noisy intermediate states.

At any timestep \(t\), the corrupted adjacency matrix \(x_t\) represents probabilistic edge states. A discretized version, \(\tilde{x}_t\), can be obtained by thresholding or rounding \(x_t\), ensuring interpretability and compatibility with downstream tasks. Both \(x_t\) and \(\tilde{x}_t\) encode the probabilities or binary states of edges, progressively transitioning from noise to structure during the reverse process.

DIFUSCO employs two distinct diffusion paradigms:
\begin{itemize}
    \item \textbf{Gaussian Diffusion (DDPM)}: Adds continuous Gaussian noise to \(x_0\), resulting in \(x_t \in \mathbb{R}^{N \times N}\).
    \item \textbf{Categorical Diffusion (D3PM)}: Uses a transition matrix to model probabilistic transitions between edge states, resulting in \(x_t \in \mathbb{R}^{N \times N \times 2}\). This approach aligns more naturally with graph-based representations and has shown superior performance in maintaining structural integrity.
\end{itemize}

\paragraph{Neural Network Architecture and Reverse Process.}  
DIFUSCO utilizes attention mechanism to guide the reverse diffusion process. The network inputs include:
\begin{itemize}
    \item The corrupted adjacency matrix \(\tilde{x}_t\),
    \item The sinusoidal encoding of the timestep \(t\),
    \item Additional node positional information, such as coordinates or learned embeddings, seamlessly integrated with edge features through an attention mechanism that effectively captures interactions between nodes and edges.
\end{itemize}

The model outputs \(\hat{x}_0 \in [0, 1]^{N \times N \times 2}\), which approximates the original adjacency matrix \(x_0\). During the reverse process, \(\hat{x}_0\) is iteratively refined and used to sample the next state \(x_{t-1}\). This iterative denoising aligns the adjacency matrix closer to the original structure \(x_0\), ensuring a gradual and consistent reconstruction.

\paragraph{Feasibility and Tour Generation.}  
DIFUSCO ensures the generation of valid Hamiltonian cycles through two main approaches: greedy decoding and sampling-based decoding. Both methods leverage the heatmap generated by the diffusion process, which contains scores \( A_{ij} \) representing the confidence of edge \((i, j)\) being part of the optimal tour. These scores are combined with node coordinates \( c_i, c_j \) to calculate a normalized ranking score:
\[
\text{Score}(i, j) = \frac{A_{ij} + A_{ji}}{\|c_i - c_j\|},
\]
where \(\|c_i - c_j\|\) is the Euclidean distance between nodes \(i\) and \(j\). The two approaches differ in how they construct the initial tour:

\subparagraph{Greedy Decoding.}  
In the greedy decoding approach, edges are sorted by their scores in descending order. The algorithm iteratively selects the highest-ranking edges while avoiding subtours and maintaining feasibility constraints, such as ensuring that each node has exactly two connections. This process ensures the formation of a valid Hamiltonian cycle. Greedy decoding employs an exhaustive 50-step inference process to refine edge probabilities before constructing the tour, prioritizing precision in edge selection.

\subparagraph{Sampling-Based Decoding.}  
In contrast, the sampling-based approach generates multiple adjacency matrices in parallel by initializing the reverse diffusion process with different random seeds. Up to 16 adjacency matrices are sampled simultaneously, with each undergoing 10 diffusion inference steps. These parallel samples provide a diverse set of candidate tours, which can then be refined using post-processing techniques to ensure feasibility and improve quality. Sampling prioritizes solution diversity, making it particularly effective in exploring a broader range of potential solutions.

Regardless of whether the initial tour is generated via greedy decoding or sampling, DIFUSCO employs advanced refinement techniques to optimize the tour further:
\begin{itemize}
    \item \textbf{2-opt Local Search:} This heuristic iteratively examines pairs of edges in the tour. If swapping two edges reduces the overall path length, the swap is performed. This process continues until no further improvements can be made, effectively eliminating inefficient connections and enhancing tour optimality.
    \item \textbf{Monte Carlo Tree Search (MCTS):} MCTS explores alternative configurations of the tour by leveraging the heatmap scores as guidance. Using \(k\)-opt transformations, it balances the exploration of new solutions and the exploitation of promising candidates. MCTS iteratively refines the tour, often leading to significant quality improvements over the initial solution.
\end{itemize}

\subsection{Diffusion Frameworks For DIFUSCO}
\subsubsection{Denoising Diffusion Probabilistic Models (DDPM)}

\paragraph{Forward Process.}  
DDPM \cite{ho2020ddpm} operates on continuous data distributions by introducing Gaussian noise in a forward process. The process starts with the original data \(x_0\), representing the final heatmap, and iteratively adds Gaussian noise at each timestep \(t\) until the data becomes random noise \(x_T\), approximated as \(x_T \sim \mathcal{N}(0, I)\). The forward process is defined as:
\[
q(x_t | x_{t-1}) = \mathcal{N}(x_t; \sqrt{1 - \beta_t} x_{t-1}, \beta_t I),
\]
where \(x_t\) represents the noisy data at timestep \(t\), and \(\beta_t \in (0, 1)\) is the variance controlling the noise level. Over multiple timesteps, the structure in \(x_0\) is progressively corrupted, resulting in \(x_T\), which serves as the starting point for the reverse process.

\paragraph{Reverse Process.}  
The reverse process reconstructs \(x_0\) from \(x_T\) by iteratively denoising through intermediate steps \(x_{T-1}, x_{T-2}, \dots, x_0\). At each timestep \(t\), the reverse distribution is parameterized as:
\[
p_\theta(x_{t-1} | x_t) = \mathcal{N}(x_{t-1}; \mu_\theta(x_t, t), \Sigma_\theta(x_t, t)),
\]
where \(\mu_\theta\) is the predicted mean, \(\Sigma_\theta = \beta_t I\) is the fixed variance, and both are predicted by a neural network. The mean \(\mu_\theta\) is computed as:
\[
\mu_\theta(x_t, t) = \frac{1}{\sqrt{\alpha_t}} \left( x_t - \frac{\beta_t}{\sqrt{1 - \alpha_t}} \epsilon_\theta(x_t, t) \right),
\]
where \(\alpha_t = \prod_{s=1}^t (1 - \beta_s)\) represents the cumulative noise schedule, and \(\epsilon_\theta(x_t, t)\) is the predicted noise at timestep \(t\). Starting from \(x_T\), the reverse process iteratively refines the noisy data, reducing the noise at each step to reconstruct \(x_0\).

The reverse process outputs either the predicted noise \(\epsilon_\theta\) at each timestep, which is used to refine \(x_t\), or the reconstructed data \(x_0\) at the end of the process. This iterative denoising ensures that \(x_0\) is a high-quality approximation of the original data distribution, making DDPM an effective generative framework.

\paragraph{Training Objective.}  
The training process optimizes the variational lower bound (VLB) to ensure that the reverse process accurately reconstructs the original data. The objective is expressed as:
\begin{align}
\mathcal{L} = \mathbb{E}_{q(x_0)} \bigg[ 
& D_{\mathrm{KL}}(q(x_T | x_0) \| p(x_T)) \nonumber \\
& + \sum_{t=2}^{T} D_{\mathrm{KL}}(q(x_{t-1} | x_t, x_0) \| p_\theta(x_{t-1} | x_t)) \nonumber \\
& - \log p_\theta(x_0 | x_1)
\bigg],
\end{align}
where \(D_{\mathrm{KL}}\) is the Kullback-Leibler divergence, \(q(x_T | x_0)\) represents the forward process’s noise distribution, and \(p_\theta(x_{t-1} | x_t)\) parameterizes the reverse process.

\paragraph{Simplified Loss.}  
For practical implementation, the training objective is often simplified to directly predict the noise \(\epsilon\) added during the forward process, leading to the following reparameterized loss:
\[
\mathcal{L}_{\mathrm{simple}} = \mathbb{E}_{t, x_0, \epsilon} \big[ \| \epsilon - \epsilon_\theta(x_t, t) \|^2 \big],
\]
where \(\epsilon \sim \mathcal{N}(0, I)\) is the ground truth noise. This simplified loss allows efficient optimization of the neural network to predict noise accurately.

\subsubsection{Discrete Denoising Diffusion Probabilistic Models (D3PM)}

D3PM \cite{austin2021structured} generalizes diffusion models to discrete state spaces, making it particularly effective for tasks where data is inherently categorical or binary. In graph-based combinatorial optimization problems, such as those addressed by DIFUSCO \cite{sun2023difusco}, the data is represented as an $N \times N \times 2$ tensor, where each entry corresponds to the state of an edge between two nodes. The tensor encodes probabilities for the two discrete states: \(k = 0\) for no edge (\(1 - p_{\text{edge}}\)) and \(k = 1\) for an edge (\(p_{\text{edge}}\)).

\paragraph{Forward Process.}  
The forward process progressively corrupts the original adjacency tensor \(x_0\) by applying a transition matrix \(Q_t \in \mathbb{R}^{2 \times 2}\) independently to each edge state at timestep \(t\). The transition matrix \(Q_t\) is defined as:
\[
Q_t = \begin{bmatrix} 
(1 - \beta_t) & \beta_t \\ 
\beta_t & (1 - \beta_t)
\end{bmatrix},
\]
where \(\beta_t \in (0, 1)\) governs the corruption level at timestep \(t\):
\begin{itemize}
    \item \((1 - \beta_t)\): Probability of maintaining the current state (edge or no edge),
    \item \(\beta_t\): Probability of transitioning to the opposite state.
\end{itemize}

At each timestep \(t\), the forward process is represented as:
\[
q(x_t | x_{t-1}) = \mathrm{Cat}(x_t; \tilde{x}_{t-1} Q_t),
\]
where \(\tilde{x}_{t-1} Q_t\) computes the probabilities of transitioning to each state based on the prior state and the corruption level \(\beta_t\). The recursion is defined as:
\[
x_t = \tilde{x}_{t-1} Q_t,
\]
applying \(Q_t\) to all entries of \(\tilde{x}_{t-1}\). Over multiple timesteps, the adjacency tensor \(x_t\) becomes increasingly corrupted, transitioning between edge and no-edge states with probabilities governed by \(\beta_t\). By the final timestep \(T\), the corrupted state \(x_T\) converges to a uniform or degenerate categorical distribution, effectively destroying all structural information in \(x_0\).

This formulation ensures systematic corruption where the degree of randomness introduced at each timestep is precisely controlled by \(\beta_t\), maintaining consistency across the \(N \times N \times 2\) tensor structure.

\paragraph{Reverse Process.}  
The reverse transition \(q(x_{t-1} | x_t, \tilde{x}_0)\) is derived using Bayes' theorem and is defined as:
\[
q(x_{t-1} | x_t, \tilde{x}_0) = \frac{q(x_t | x_{t-1}, \tilde{x}_0) q(x_{t-1} | \tilde{x}_0)}{q(x_t | \tilde{x}_0)},
\]
where \(q(x_t | x_{t-1}, \tilde{x}_0)\) represents the forward corruption step from \(x_{t-1}\) to \(x_t\), \(q(x_{t-1} | \tilde{x}_0)\) encodes the marginal transition from \(\tilde{x}_0)\) to \(x_{t-1}\) and \(q(x_t | \tilde{x}_0)\) normalizes the posterior over \(x_{t-1}\).

The forward transition across multiple timesteps is captured by the cumulative transition matrix \(\bar{Q}_t = Q_1 Q_2 \cdots Q_t\), where each \(Q_t \in \mathbb{R}^{2 \times 2}\) represents the corruption probabilities at timestep \(t\). The entries of \(\bar{Q}_t\) encode the overall probabilities of transitioning between the two states (edge or no edge) over the first \(t\) steps. This allows the marginal \(q(x_t | \tilde{x}_0)\) to be expressed as:
\[
q(x_t | \tilde{x}_0) = \mathrm{Cat}\left( x_t; p = \tilde{x}_0 \bar{Q}_t \right).
\]

Using this cumulative formulation, the posterior distribution \(q(x_{t-1} | x_t, \tilde{x}_0)\) is given by:
\[
q(x_{t-1} | x_t, \tilde{x}_0) = \mathrm{Cat}\left( x_{t-1}; p = \frac{\tilde{x}_t Q_t^\top \odot \tilde{x}_0 \bar{Q}_{t-1}}{\tilde{x}_0 \bar{Q}_t \tilde{x}_t^\top} \right),
\]
where \(\tilde{x}_t\) is the discretized adjacency tensor at timestep \(t\), \(Q_t^\top\) is the transpose of the transition matrix for timestep \(t\) and \(\odot\) denotes element-wise multiplication.

This formulation ensures that the reverse sampling process systematically reconstructs \(x_0\) while incorporating both the corrupted state \(x_t\) and the clean prediction \(\tilde{x}_0\).

The neural network \(p_\theta\) is trained to predict \(\tilde{x}_{0, \theta}\) (discretized version of \(x_{0, \theta}\)), the clean adjacency tensor, given the corrupted \(x_t\). The reverse process then substitutes the predicted \(\tilde{x}_0\) into the posterior to compute:
\[
p_\theta(x_{t-1} | x_t) = \sum_{\tilde{x}_0} q(x_{t-1} | x_t, \tilde{x}_0) p_\theta(\tilde{x}_0 | x_t).
\]

By iteratively applying this reverse sampling from \(x_T\) to \(x_0\), the process reconstructs the original adjacency tensor while maintaining consistency with the forward corruption model. The integration of the cumulative transition matrix \(\bar{Q}_t\) ensures that the overall corruption and reconstruction remain aligned across all timesteps.

\paragraph{Training Objective.}  
The model is trained to minimize the binary cross-entropy loss between the predicted discretized tensor \(\hat{x}_0\) (output by the neural network) and the ground truth discretized tensor \(x_0\). The process involves uniformly sampling a timestep \(t \in \{1, 2, \dots, T\}\), applying the forward corruption process to obtain \(x_t\), and using the neural network to predict \(x_{0, \theta} (x_t)\) from \(x_t\). The loss function is defined as:
\begin{align}
\mathcal{L} = \mathbb{E}_{t, x_0} \bigg[ 
& -\sum_{i=1}^N \sum_{j=1}^N \sum_{k=0}^1 \big( 
x_0(i, j, k) \log x_{0, \theta}(i, j, k) \nonumber \\
& + (1 - x_0(i, j, k)) \log (1 - x_{0, \theta}(i, j, k)) 
\big) 
\bigg],
\end{align}
where \(x_0(i, j, k) \in \{0, 1\}\) is the ground truth binary state for the edge between nodes \(i\) and \(j\) and \(x_{0, \theta}(i, j, k)\) is the predicted probability output by the neural network for the same edge state.

The training process ensures that the neural network learns to predict the ground truth discretized tensor \(x_0\) by minimizing the binary cross-entropy loss. By comparing \(x_{0, \theta}\) and \(x_0\) directly, the model aligns its predictions with the underlying structure of the original adjacency tensor, enabling accurate reconstruction during the reverse process.

\section{Blackout DIFUSCO: Proposition}
\subsection{Motivation.}  
In DIFUSCO, the discrete diffusion model (D3PM) outperformed its continuous Gaussian counterpart (DDPM) due to its ability to naturally align with the discrete nature of adjacency matrix. Inspired by this success, we propose integrating Blackout Diffusion into DIFUSCO. Blackout Diffusion extends the discrete state-space paradigm to a continuous-time framework, combining the structural advantages of D3PM with the smoother transitions and finer control offered by continuous-time dynamics. This integration is anticipated to further enhance the model's ability to generate high-quality solutions for combinatorial optimization problems like the Traveling Salesman Problem (TSP).

\subsection{Blackout DIFUSCO}
\paragraph{Forward Process.}  
The forward process in Blackout Diffusion progressively corrupts the original adjacency matrix \( x_0 \in \{0, 1\}^{N \times N} \), representing the binary edge connectivity, into a fully degenerate blackout state \( x_T = 0 \). The dynamics are governed by a pure-death Markov process, where each edge state transitions from active (\(1\), edge is in the tour) to inactive (\(0\), edge is not in the tour) over a continuous time interval \([0, T]\). The forward process can be expressed as:
\[
q(x_t | x_0) = \mathrm{Binomial}(x_t; x_0, e^{-t}),
\]
where \( e^{-t} \) represents the decay rate of edge states over time. 

As \( t \to T \), the corrupted matrix \( x_t \) converges to a black image, effectively destroying all structural information in \( x_0 \). To ensure complete corruption, we set \( T \) to a sufficiently large value of 15.0, as done in this work.

\paragraph{Reverse Process.}  
The reverse process reconstructs the adjacency matrix \( x_0 \) by simulating a birth-only process, effectively reversing the pure-death dynamics of the forward process. Each element of the matrix, corresponding to an edge state \((i, j)\), is processed independently. The reverse transition probability for edge state \((i, j)\) is modeled as:

\begin{equation}
p(x_{s}(i, j) | x_t(i, j), x_0(i, j)) = 
\binom{o - n}{m - n} \cdot r^{m-n} (1 - r)^{o-m},
\label{eq:blackout_reverse}
\end{equation}

where the variables are defined as follows:
\begin{itemize}
    \item \( n = x_t(i, j) \): The corrupted edge state at timestep \( t \),
    \item \( m = x_{s}(i, j) \): The intermediate reconstructed edge state at timestep \( s \),
    \item \( o = x_0(i, j) \): The original binary edge state,
    \item \( r = \frac{e^{-s} - e^{-t}}{1 - e^{-t}} \): The binomial transition parameter.
\end{itemize}

This formulation, derived from the \textit{binomial bridge}, ensures a smooth and probabilistically consistent transition from the corrupted state \( x_t \) to the original state \( x_0 \).

To approximate \( \hat{x}_0 \) from \( \hat{x}_t \), a neural network \(\text{NN}_\theta\) is used to predict the change in state \(\Delta x_{t \to 0}\), representing the difference between \( \hat{x}_t \) and the original state \( \hat{x}_0 \), as:
\[
\Delta x_{t \to 0} = \text{NN}_\theta(\hat{x}_t, t)
\]
where:
\begin{itemize}
    \item \( \text{NN}_\theta \): A parameterized neural network,
    \item Inputs: The corrupted adjacency matrix \( \hat{x}_t \) and the timestep \( t \),
    \item Outputs: The predicted change \(\Delta x_{t \to 0}\) between \( \hat{x}_t \) and \( \hat{x}_0 \).
\end{itemize}

The reconstruction of \( x_s \) is then performed iteratively as:
\[
\hat{x}_0 = \hat{x}_t + \Delta x_{t \to 0}.
\]
The probability \( p(x_{s} | \hat{x}_t, \hat{x}_0) \) follows a similar formulation to the binomial transition probability shown in Equation~\ref{eq:blackout_reverse}, ensuring consistency with the dynamics of the reverse process.

This iterative refinement ensures that the adjacency matrix \( x_0 \) is reconstructed accurately, element by element, while respecting the probabilistic structure modeled by the binomial bridge. Each edge state evolves independently during the reverse process, guided by the neural network's predictions and the binomial transition probabilities.

\paragraph{Training Objective.}  
The training objective is designed to minimize the discrepancy between the predicted and true adjacency matrice, incorporating the temporal dynamics of the diffusion process. The loss function for each element \(i\) is formulated as:
\begin{equation}
l_i = (t_k - t_{k-1}) e^{-t_k} \left[ y_i - \left( (x_0 - x_{t_k})_i \log y_i \right) \right],
\label{eq:training_objective}
\end{equation}
where:
\begin{itemize}
    \item \(t_k\) and \(t_{k-1}\) represent successive timesteps,
    \item \(y_i\) is the difference between two states for the element \(i\), which is \(\Delta x_{t \to 0}\),
    \item \(x_0\) is the ground truth adjacency matrix,
    \item \(x_{t_k}\) is the corrupted adjacency matrix at timestep \(t_k\).
\end{itemize}

This loss function ensures that the model learns to predict transitions in the adjacency matrix while respecting the temporal dynamics introduced during the diffusion process. By incorporating the time-dependent weight \(e^{-t_k}\) and the logarithmic term \(\log y_i\), the objective penalizes large deviations in the reconstructed probabilities and encourages convergence towards the true solution.

\paragraph{Simplified Loss.}  
The training objective can also be simplified to focus directly on the discrepancy between the predicted and true adjacency tensors, with less emphasis on the temporal dynamics. The simplified loss function for each element \(i\) is expressed as:
\begin{equation}
l_i = \left[ y_i - \left( (x_0 - x_{t_k})_i \log y_i \right) \right],
\label{eq:Simplified_Loss}
\end{equation}

This formulation simplifies the original loss by removing the explicit time-dependent weighting terms, allowing the focus to remain on accurately reconstructing the true adjacency tensor while still leveraging the diffusion dynamics for guidance.

\begin{figure*}[ht]
    \centering
    \begin{minipage}{0.32\textwidth}
        \centering
        \includegraphics[width=\linewidth]{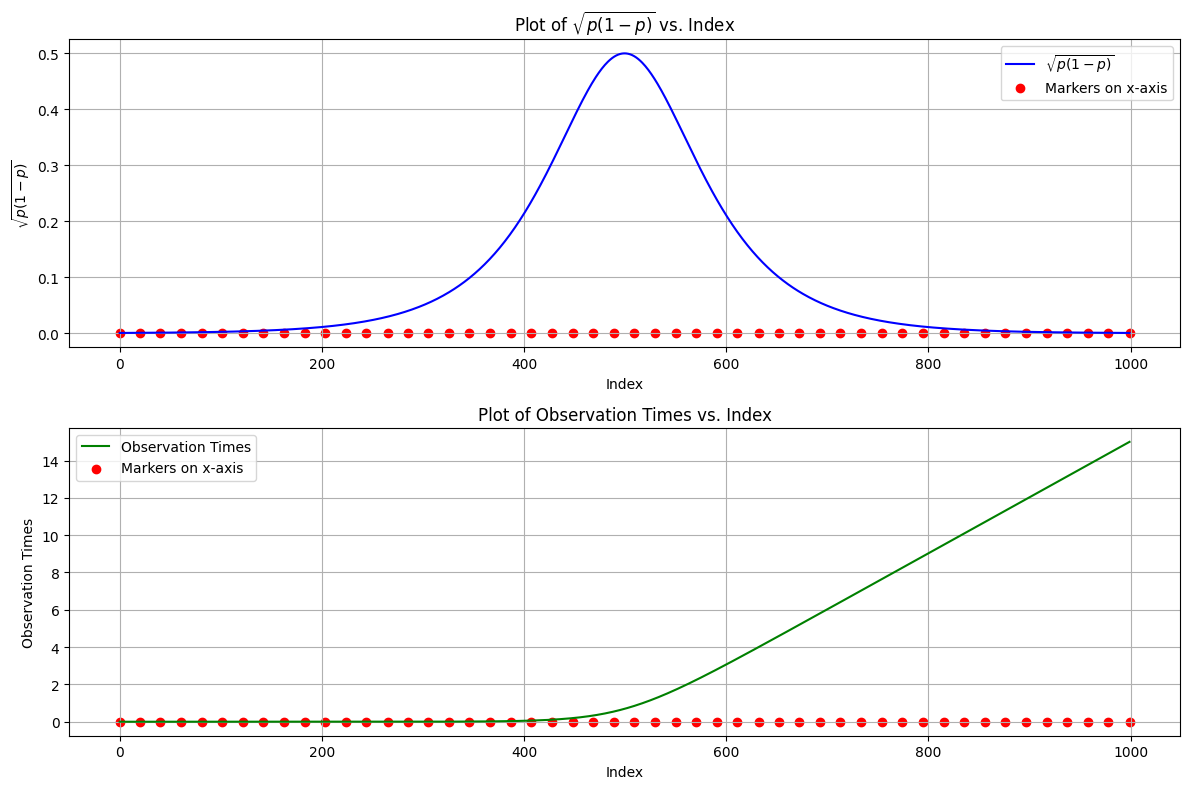}
        \caption*{(a) Original Blackout DIFUSCO}
    \end{minipage}
    \begin{minipage}{0.32\textwidth}
        \centering
        \includegraphics[width=\linewidth]{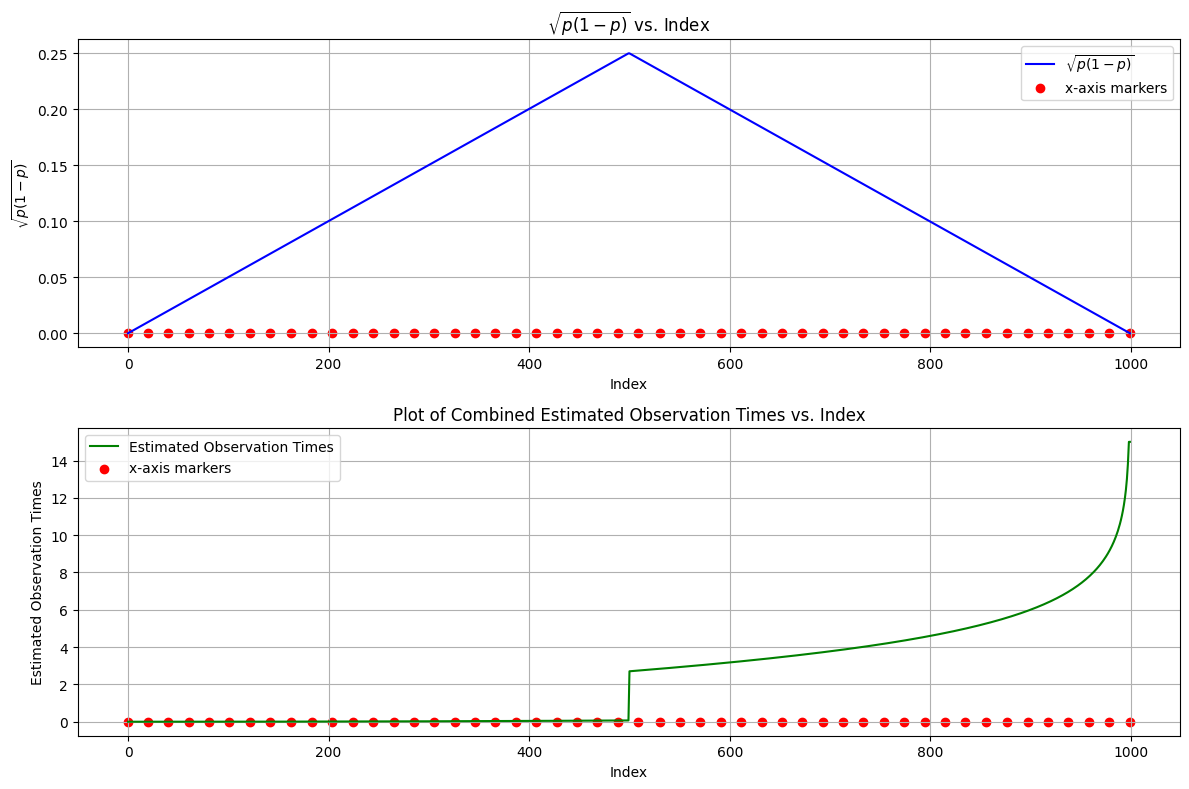}
        \caption*{(b) Improved Blackout DIFUSCO}
    \end{minipage}
    \begin{minipage}{0.32\textwidth}
        \centering
        \includegraphics[width=\linewidth]{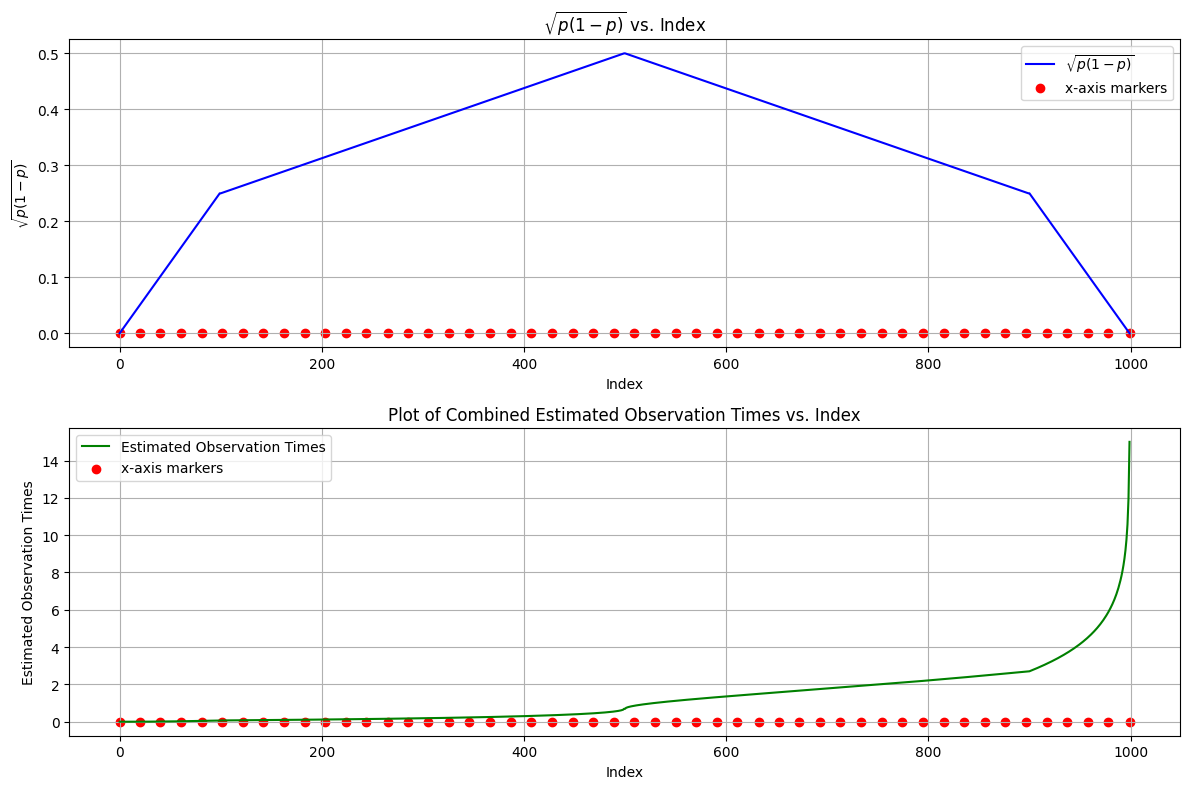}
        \caption*{(c) More Improved Blackout DIFUSCO}
    \end{minipage}
    \caption{
        Visualization of the forward diffusion process in Blackout DIFUSCO variants. 
        The plots illustrate the evolution of the corrupted heatmaps as the process progresses across the time axis. 
        Each design leverages the binomial distribution for modeling the forward process, and the variance of the heatmaps (\(\text{std}\)) is analyzed to optimize observation time scheduling. The red markers in the plots indicate the sampling points during the generation step, highlighting the regions where the model focuses its computations.
        \textbf{(a) Original Blackout DIFUSCO}: Displays the standard forward process with no modifications.
        \textbf{(b) Improved Blackout DIFUSCO}: Implements a peaked observation time schedule inspired by the cosine scheduler, ensuring that the standard deviation grows linearly to a sharp midpoint and tapers symmetrically.
        \textbf{(c) More Improved Blackout DIFUSCO}: Further enhances the scheduling by introducing a hyperparameter \(\alpha = 0.2\), which redistributes sampling density to focus on regions with maximum \(\text{std}\), ensuring improved reconstruction of challenging adjacency matrix states.
    }
    \label{fig:blackout_difusco_comparison}
\end{figure*}

\subsection{Improved Blackout DIFUSCO}
To enhance the performance of Blackout DIFUSCO, we propose a modification to the observation time scheduling inspired by the cosine scheduler from Improved DDPM \cite{nichol2021improved}. Specifically, we analyze the standard deviation (\(\text{std}\)) of the corrupted heatmaps across the time axis and observe that, unlike Gaussian-based diffusion models, the variance in Blackout Diffusion is smallest at the final corruption stage (\(t = T\)) and largest at the midpoint. 

To address this, we redesign the observation time scheduling such that the standard deviation increases linearly throughout the diffusion process, producing a sharper peak at intermediate timesteps. This adjustment ensures that the model allocates more computational focus during the middle stages, where the variance is highest, and less attention towards the extremes, where the variance is minimal. This "peaked" observation time distribution aligns with the unique characteristics of Blackout Diffusion, allowing the reverse process to reconstruct adjacency matrices more effectively.

\subsection{More Improved Blackout DIFUSCO}
Building upon the observation time scheduling introduced in Improved Blackout DIFUSCO, we further refine the design to focus on regions where the corrupted heatmaps exhibit the largest standard deviations (\(\text{std}\)), as these regions are more challenging for the network to reconstruct. To achieve this, we introduce a hyperparameter \(\alpha\) to control the sampling density along the time axis.

The new schedule is designed to allocate more sampling steps to regions with high \(\text{std}\) while maintaining symmetry around the point of maximum variance. Specifically:
\begin{itemize}
    \item The middle portion of the time axis, where the variance is highest, receives more samples.
    \item The start and end segments (\(T \alpha\)) are linearly reduced to half the maximum variance, creating a smooth transition.
\end{itemize}

For our experiments, we set \(\alpha = 0.2\), ensuring that \(20\%\) of the time axis on both ends tapers to half the maximum \(\text{std}\). The middle portion transitions linearly from this midpoint back to the peak, forming a symmetric triangular schedule. This design enables the model to concentrate more sampling steps in regions with high reconstruction difficulty, further improving its ability to recover accurate adjacency matrices.

\section{Experiments}

\subsection{Dataset Preparation}
The models were trained on the TSP-50 dataset, consisting of 1,502,000 training samples, and evaluated on a test set of 1,280 samples. To assess the generalization capabilities, the trained models were further evaluated on TSP-100 and TSP-500 datasets, with test set sizes of 1,280 and 128 samples, respectively. All experiments were conducted on a system equipped with an NVIDIA GeForce RTX 3090 GPU.

\subsection{Comparison with Other Models}

\begin{table}[ht]
    \centering
    \caption{Comparison of Greedy-Based Approaches with Traditional Solvers on the TSP-50 Dataset. Concorde and 2-OPT are traditional solvers, with Concorde being capable of providing exact solutions to the TSP. The other methods are evaluated based on their greedy approaches. The Top-1, Top-2, and Top-3 gaps are underlined for clarity.}
    \label{tab:comparison}
    \begin{tabular}{|l|c|c|}
        \hline  
        \textbf{Algorithm} & \textbf{Length} & \textbf{Gap (\%)} \\ \hline
        Concorde              & 5.69   & 0.00   \\
        2-OPT                  & 5.86   & 2.95   \\
        \hline
        AM \cite{kool2019attention}                     & 5.80   & 1.76   \\
        GCN \cite{joshi2019tsp}                   & 5.87   & 3.10   \\
        Transformer \cite{bresson2021transformer}       & 5.71   & 0.31   \\
        POMO \cite{kwon2020pomo}                   & 5.73   & 0.64   \\
        Sym-NCO \cite{kim2022symnco}                & 5.74   & 0.88   \\
        Image Diffusion \cite{song2021ddim}       & 5.76   & 1.23   \\
        DIFUSCO (Greedy) \cite{sun2023difusco}          & 5.70   &  \underline{\textbf{0.17}}   \\
        T2T (Greedy) \cite{li2023t2t}        & 5.70   & \underline{\textbf{0.10}}   \\ \hline
        B-DIFUSCO (original) & 5.70 & 0.23  \\
        B-DIFUSCO (improved) & 5.71 & 0.38  \\
        B-DIFUSCO (more improved) & 5.70 & \underline{\textbf{0.20}}  \\
        \hline
    \end{tabular}
\end{table}

Table~\ref{tab:comparison} outlines the performance of the Blackout DIFUSCO variants and other greedy-based approaches alongside traditional solvers like Concorde and 2-OPT. The results highlight:
\begin{itemize}
    \item The \textbf{Top-3 Greedy Performers}: DIFUSCO (Greedy), T2T (Greedy), and B-DIFUSCO (more improved), which demonstrated the lowest gaps (\%) compared to other methods.
    \item Blackout DIFUSCO's performance placed it consistently within the \textbf{Top-3}, emphasizing its competitiveness and efficiency compared to the baseline models.
\end{itemize}

\subsection{Baseline Comparisons}

\begin{figure}[ht]
    \centering
    \includegraphics[width=0.5\textwidth]{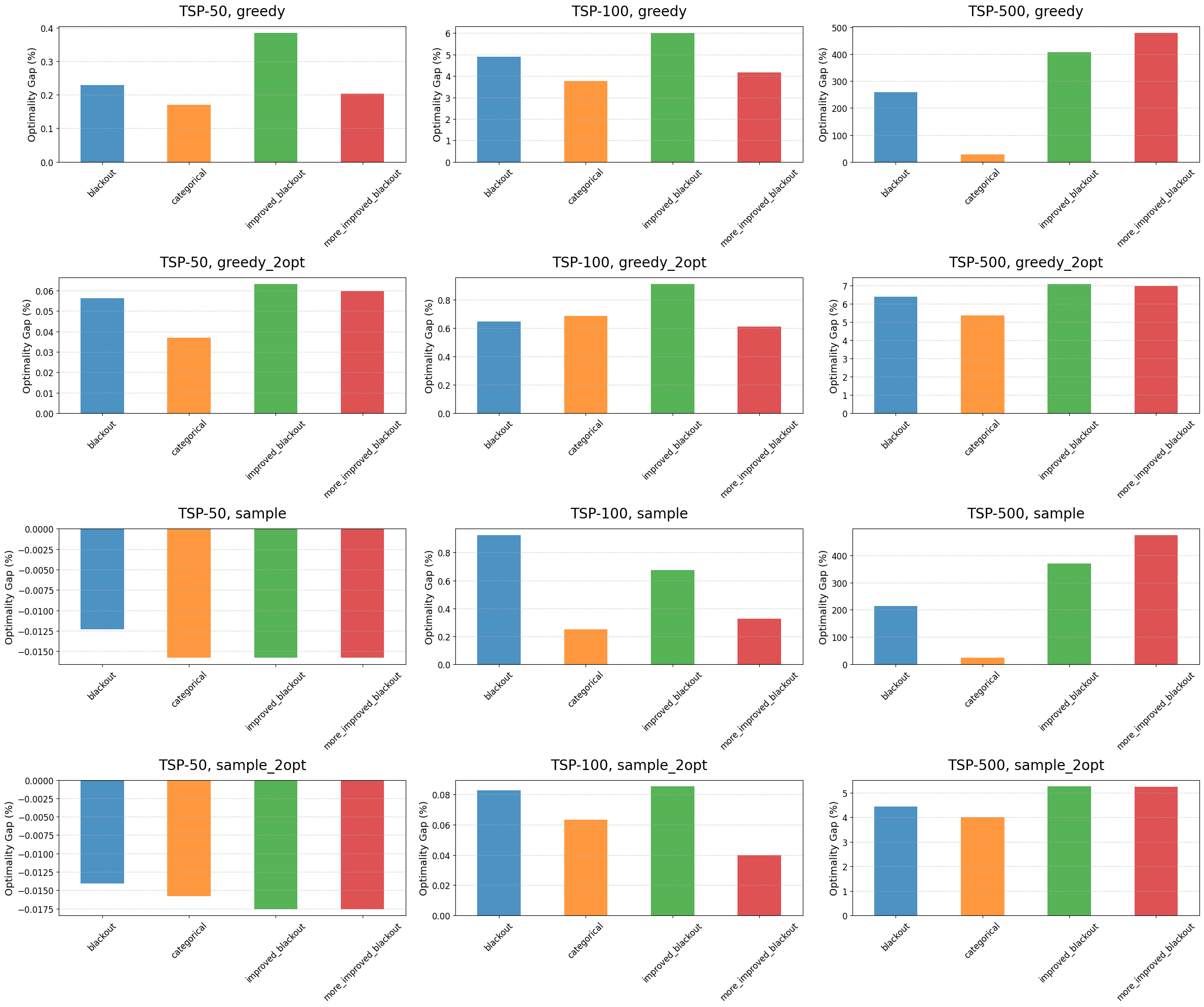}
    \caption{Bar plots showing Optimality Gap (\%) for different decoding methods (Greedy, Greedy + 2-OPT, Sampling, Sampling + 2-OPT) across datasets (TSP-50, TSP-100, TSP-500) and categories (Categorical DIFUSCO, Blackout DIFUSCO, Improved Blackout DIFUSCO, More Improved Blackout DIFUSCO).}
    \label{fig:optimality_gap}
\end{figure}

Table~\ref{tab:comparison} also demonstrates the progression in performance among Blackout DIFUSCO variants:
\begin{itemize}
    \item \textbf{Original Blackout DIFUSCO} achieved a gap of \textbf{0.23\%}, indicating a strong baseline.
    \item \textbf{Improved Blackout DIFUSCO} further enhanced the performance but slightly increased the gap to \textbf{0.38\%} due to additional complexity.
    \item \textbf{More Improved Blackout DIFUSCO} finalized with a gap of \textbf{0.20\%}, achieving near-optimal performance within the greedy-based category.
\end{itemize}
This incremental improvement underscores the effectiveness of targeted refinements in the diffusion model.

\subsection{Result Analysis}
Figure~\ref{fig:optimality_gap} provides a visual analysis of the optimality gaps across decoding methods (Greedy, Greedy + 2-OPT, Sampling, Sampling + 2-OPT) and datasets (TSP-50, TSP-100, TSP-500). Notable insights include:
\begin{itemize}
    \item \textbf{Generalization Capability}: Although trained exclusively on TSP-50, the Blackout DIFUSCO model effectively generalized to larger datasets (TSP-100 and TSP-500), maintaining competitive gaps.
    \item \textbf{Model Stability}: The consistent performance of the diffusion-based methods across all datasets suggests robustness against instance complexity.
    \item \textbf{Greedy Efficiency}: The results reaffirm the strength of greedy strategies when optimized diffusion techniques are applied.
\end{itemize}

Overall, the results validate the proposed methods' capabilities and their potential for solving combinatorial optimization problems efficiently.

\section{Conclusion}
The integration of Blackout Diffusion into DIFUSCO introduces a novel approach to solving combinatorial optimization problems through continuous-time dynamics. By aligning variance peaks with observation time and concentrating sampling in regions of high reconstruction difficulty, the proposed improvements enhance solution quality and model efficiency. While More Improved Blackout Diffusion demonstrated superior performance on TSP50 and competitive results on larger datasets, the method underscores the potential of continuous-time frameworks in combinatorial optimization. This study serves as a foundational step towards advancing diffusion-based models for discrete optimization tasks, encouraging future exploration of optimized scheduling techniques and adaptive sampling strategies. The insights gained here pave the way for broader applications of generative diffusion models in solving NP-hard problems.
\label{sec:rationale}

\clearpage

\setcounter{page}{1}
\section*{Supplementary Material}

\subsection*{TSP 50 - Metrics by Post-Processing Method}

\begin{table}[H]
    \centering
    \caption{TSP Metrics for Original Blackout DIFUSCO (TSP 50)}
    \begin{tabular}{lrr}
        \toprule
        Post-Processing Method & val\_solved\_cost & Opt Gap (\%) \\
        \midrule
        GREEDY                 & 5.7018 & 0.23 \\
        GREEDY + 2-OPT         & 5.6920 & 0.06 \\
        SAMPLE                 & 5.6881 & -0.01 \\
        SAMPLE + 2-OPT         & 5.6879 & -0.02 \\
        \bottomrule
    \end{tabular}
\end{table}

\begin{table}[H]
    \centering
    \caption{TSP Metrics for Improved Blackout DIFUSCO (TSP 50)}
    \begin{tabular}{lrr}
        \toprule
        Post-Processing Method & val\_solved\_cost & Opt Gap (\%) \\
        \midrule
        GREEDY                 & 5.7107 & 0.38 \\
        GREEDY + 2-OPT         & 5.6924 & 0.06 \\
        SAMPLE                 & 5.6879 & -0.02 \\
        SAMPLE + 2-OPT         & 5.6878 & -0.02 \\
        \bottomrule
    \end{tabular}
\end{table}

\begin{table}[H]
    \centering
    \caption{TSP Metrics for More Improved Blackout DIFUSCO (TSP 50)}
    \begin{tabular}{lrr}
        \toprule
        Post-Processing Method & val\_solved\_cost & Opt Gap (\%) \\
        \midrule
        GREEDY                 & 5.7004 & 0.20 \\
        GREEDY + 2-OPT         & 5.6922 & 0.06 \\
        SAMPLE                 & 5.6879 & -0.02 \\
        SAMPLE + 2-OPT         & 5.6878 & -0.02 \\
        \bottomrule
    \end{tabular}
\end{table}

\begin{table}[H]
    \centering
    \caption{TSP Metrics for Categorical DIFUSCO (TSP 50)}
    \begin{tabular}{lrr}
        \toprule
        Post-Processing Method & val\_solved\_cost & Opt Gap (\%) \\
        \midrule
        GREEDY                 & 5.6985 & 0.17 \\
        GREEDY + 2-OPT         & 5.6909 & 0.04 \\
        SAMPLE                 & 5.6879 & -0.02 \\
        SAMPLE + 2-OPT         & 5.6879 & -0.02 \\
        \bottomrule
    \end{tabular}
\end{table}

\subsection*{TSP 100 - Metrics by Post-Processing Method}

\begin{table}[H]
    \centering
    \caption{TSP Metrics for Original Blackout DIFUSCO (TSP 100)}
    \begin{tabular}{lrr}
        \toprule
        Post-Processing Method & val\_solved\_cost & Opt Gap (\%) \\
        \midrule
        GREEDY                 & 8.1402 & 4.91 \\
        GREEDY + 2-OPT         & 7.8095 & 0.65 \\
        SAMPLE                 & 7.8310 & 0.92 \\
        SAMPLE + 2-OPT         & 7.7657 & 0.08 \\
        \bottomrule
    \end{tabular}
\end{table}

\begin{table}[H]
    \centering
    \caption{TSP Metrics for Improved Blackout DIFUSCO (TSP 100)}
    \begin{tabular}{lrr}
        \toprule
        Post-Processing Method & val\_solved\_cost & Opt Gap (\%) \\
        \midrule
        GREEDY                 & 8.2264 & 6.02 \\
        GREEDY + 2-OPT         & 7.8301 & 0.91 \\
        SAMPLE                 & 7.8115 & 0.67 \\
        SAMPLE + 2-OPT         & 7.7659 & 0.08 \\
        \bottomrule
    \end{tabular}
\end{table}

\begin{table}[H]
    \centering
    \caption{TSP Metrics for More Improved Blackout DIFUSCO (TSP 100)}
    \begin{tabular}{lrr}
        \toprule
        Post-Processing Method & val\_solved\_cost & Opt Gap (\%) \\
        \midrule
        GREEDY                 & 8.0820 & 4.16 \\
        GREEDY + 2-OPT         & 7.8067 & 0.61 \\
        SAMPLE                 & 7.7846 & 0.33 \\
        SAMPLE + 2-OPT         & 7.7624 & 0.04 \\
        \bottomrule
    \end{tabular}
\end{table}

\begin{table}[H]
    \centering
    \caption{TSP Metrics for Categorical DIFUSCO (TSP 100)}
    \begin{tabular}{lrr}
        \toprule
        Post-Processing Method & val\_solved\_cost & Opt Gap (\%) \\
        \midrule
        GREEDY                 & 8.0516 & 3.77 \\
        GREEDY + 2-OPT         & 7.8125 & 0.69 \\
        SAMPLE                 & 7.7787 & 0.25 \\
        SAMPLE + 2-OPT         & 7.7642 & 0.06 \\
        \bottomrule
    \end{tabular}
\end{table}

\subsection*{TSP 500 - Metrics by Post-Processing Method}

\begin{table}[H]
    \centering
    \caption{TSP Metrics for Original Blackout DIFUSCO (TSP 500)}
    \begin{tabular}{lrr}
        \toprule
        Post-Processing Method & val\_solved\_cost & Opt Gap (\%) \\
        \midrule
        GREEDY                 & 59.5824 & 259.29 \\
        GREEDY + 2-OPT         & 17.6443 & 6.4 \\
        SAMPLE                 & 51.9534 & 213.28 \\
        SAMPLE + 2-OPT         & 17.3200 & 4.44 \\
        \bottomrule
    \end{tabular}
\end{table}

\begin{table}[H]
    \centering
    \caption{TSP Metrics for Improved Blackout DIFUSCO (TSP 500)}
    \begin{tabular}{lrr}
        \toprule
        Post-Processing Method & val\_solved\_cost & Opt Gap (\%) \\
        \midrule
        GREEDY                 & 84.0509 & 406.83 \\
        GREEDY + 2-OPT         & 17.7599 & 7.09 \\
        SAMPLE                 & 77.8347 & 369.35 \\
        SAMPLE + 2-OPT         & 17.4551 & 5.26 \\
        \bottomrule
    \end{tabular}
\end{table}

\begin{table}[H]
    \centering
    \caption{TSP Metrics for More Improved Blackout DIFUSCO (TSP 500)}
    \begin{tabular}{lrr}
        \toprule
        Post-Processing Method & val\_solved\_cost & Opt Gap (\%) \\
        \midrule
        GREEDY                 & 96.1816 & 479.98 \\
        GREEDY + 2-OPT         & 17.7413 & 6.98 \\
        SAMPLE                 & 95.0985 & 473.45 \\
        SAMPLE + 2-OPT         & 17.4543 & 5.25 \\
        \bottomrule
    \end{tabular}
\end{table}

\begin{table}[H]
    \centering
    \caption{TSP Metrics for Categorical DIFUSCO (TSP 500)}
    \begin{tabular}{lrr}
        \toprule
        Post-Processing Method & val\_solved\_cost & Opt Gap (\%) \\
        \midrule
        GREEDY                 & 21.3278 & 28.61 \\
        GREEDY + 2-OPT         & 17.4737 & 5.37 \\
        SAMPLE                 & 20.5145 & 23.7 \\
        SAMPLE + 2-OPT         & 17.2468 & 4.0 \\
        \bottomrule
    \end{tabular}
\end{table}

\end{document}